\documentclass[11pt]{article}
\usepackage{amsfonts}
\usepackage{bbm}
\usepackage{mathrsfs}
\usepackage{amssymb}
\usepackage{url}
\title{A Note on Regularities of The Standard Quantized Matrix Algebra $M_q(n)$\thanks{Project supported by
the Fundamental Research Funds for Universities
of Xinjiang Uygur Autonomous Region.}}
\author{Rabigul Tuniyaz\\
{\small Department of Mathematics, School of  Science}\\
{\small Xinjiang Institute of Science and Technology}\\
{\small  Akesu, 843100, Xinjiang, China}\\}

\date{}

\begin{document}
\maketitle

\begin{center}
\begin{minipage}{120mm}
{\small {\bf Abstract.} Let $M_q(n)$ be the standard quantized matrix algebra (introduced by Faddeev, Reshetikhin, and Takhtajan).  It is shown that $M_q(n)$ is Auslander regular, Cohen-Macaulay,  Artin-Schelter regular, and a maximal order in its quotient division algebra.}
\end{minipage}\end{center} {\parindent=0pt\par

}
\vskip .5truecm
\renewcommand{\thefootnote}{\fnsymbol{footnote}}
\let\footnote\relax\footnotetext{E-mail: rabigul802@sina.com}
\let\footnote\relax
\def\NZ{\mathbb{N}}
\def\QED{\hfill{$\Box$}}
\def \r{\rightarrow}
\def\mapright#1#2{\smash{\mathop{\longrightarrow}\limits^{#1}_{#2}}}
\def\v5{\vskip .5truecm}
\def\OV#1{\overline {#1}}
\def\hang{\hangindent\parindent}
\def\textindent#1{\indent\llap{#1\enspace}\ignorespaces}
\def\item{\par\hang\textindent}
\def\LH{{\bf LH}}\def\LM{{\bf LM}}\def\LT{{\bf
LT}}\def\KX{K\langle X\rangle} \def\KZ{K\langle Z\rangle}
\def\B{{\cal B}} \def\LC{{\bf LC}} \def\G{{\cal G}}
\def\PRC{\prec_{d\textrm{\tiny -}rlex}}

\section*{1. Introduction and preliminaries}
It is well known that the standard quantized matrix algebra $M_q(n)$ introduced by Faddeev, Reshetikhin, and Takhtajan in [2] has played a very important role in the theory of quantum groups and has been extensively investigated from the  perspectives  of representation theory and noncommutative geometry. Recently, some structural properties of the quantized matrix algebra $M_q(n)$ and its modules have been established in a constructive-computational way in [10]. More precisely, it was shown explicitly that the defining relations of $M_q(n)$ constitute a Gr\"obner-Shirshov basis, and that $M_q(n)$ is a solvable polynomial algebra in the sense of [3]. Consequently, a comprehensive investigation showed that $M_q(n)$ is a Noetherian domain and classical quadratic Koszul algebra with Hilbert series $\frac{1}{(1-t)^{n^2}}$, GK dimension $n^2$, and global homological dimension $n^2$, and that finite free resolutions of finitely generated modules over $M_q(n)$ can be constructed algorithmically,  the stability of finitely generated projective modules and the $K_0$-groups of $M_q(n)$ can be established, minimal graded generating sets of finitely generated graded modules may be computed algorithmically, and $M_q(n)$ has the elimination property for one-sided ideals and finitely generated modules  in the sense of ([6], [7]).\par

Based on [10], in this note we aim to show further that $M_q(n)$ is Auslander regular, Cohen-Macaulay,  Artin-Schelter regular, and a maximal order in its division quotient algebra (Theorem 3.1). To this end, we need some preliminaries as mentioned below.\v5

Throughout this note, $K$ denotes a field of characteristic 0, $K^*=K-\{0\}$, and all rings will be $K$-algebras with multiplicative identity 1. If $S$ is a nonempty subset of an algebra $A$, then we write $\langle S\rangle$ for the two-sided ideal of $A$ generated by $S$.\v5

Recall from the literature (for example, [8], [4], [1], [9]) the following\v5

{\bf Definition 1.1} (1) A $K$-algebra $A$ is called {\it Auslander-regular} provided that $A$ has finite global homological dimension and satisfies the Gorenstein condition: $j(N)\ge j$ for every finitely generated left $A$-module $M$, every integer $j\ge 0$ and every (right) $A$-submodule $N$ of Ext$^j_A(M,A)$, where $j(N)$ is the grade number of $N$ which is the least integer $i$ such that Ext$^i_A(M,A)\ne 0$.\par
(2) An Auslander-regular $K$-algebra $A$ is called {\it Cohen-Macaulay} (or CM as abbreviation)  provided that $A$ has  finite Gelfand-Kirillov dimension $n$ and,  for every finitely generated left $A$-module $M$, the equality GK.dim$M+j(M)=n$ holds, where GK.dim denotes the Gelfand-Kirollov dimension of $M$.\par
(3) Let $A$ be a (connected) $\NZ$-graded $K$-algebra,
generated in degree one, i.e., $A=\oplus_{i\in\NZ}A_i$, where $A_0=K$ is central, dim$A_i<\infty$ for all $i$, and A is generated as an algebra by $A_1$. Let $A^+=\oplus_{i\ge 1}A_i$ be the graded
radical of $A$. If $A$ has global homological dimension gl.dim$A=d$ and finite Gelfand-Kirillov dimension, such that ${\rm Ext}^n_A(K_A,A)=\delta_{nd}(_AK)$, where $K$ is identified with $A/A^+$, then $A$ is called an Artin-Schelter regular algebra of dimension $d$.\par

(4) A $K$-algebra $A$ is called stably free if, for every
finitely generated projective $A$-module $P$, there exist integers $n$ and $m$ such that
$P\oplus A^n=A^m$.\par

(5) Let $A$ be a domain and $Q(A)$ be its quotient division ring. If $A\subseteq S$ for some $K$-algebra $S$ with the property that $aSb\subseteq A$ for some $a,b\in A-\{ 0\}$ implies  $A=S$,  then $A$ is called a {\it maximal order} in  $Q(A)$.\v5

{\bf Proposition 1.2} [9, Theorem 2.10]  Let $A$ be a Noetherian K-algebra that is Auslander regular and CM. Assume that $A$ is also stably free. Then, $A$ is a domain and a maximal order in its quotient division ring $Q(A)$. \v5

Recall from [2] that the standard quantized matrix algebra $M_q(n)$ is the associative $K$-algebra generated by $n^2$ generators $z_{ij}$, where $(i,j)\in I(n)=\{(i,j)~|~i,j=1,2,\cdots,n\}$ and $n\ge 2$,  subject to the relations:
$$\begin{array}{rll}
R_{[j<k]}:&z_{ij}z_{ik}=qz_{ik}z_{ij},&\hbox{if}~j<k ,\\
R_{[i<k]}:&z_{ij}z_{kj}=qz_{kj}z_{ij},&\hbox{if}~i<k ,\\
R_{[i<s,t<j]}:&z_{ij}z_{st}=z_{st}z_{ij},&\hbox{if}~i<s, t<j,\\
R_{[i<s,j<t]}:&z_{ij}z_{st}=z_{st}z_{ij}+(q-q^{-1})z_{it}z_{sj},&\hbox{if}~i<s,~j<t,
 \end{array}$$
where $i,j,k,s,t=1,2,...,n$ and $q\in K^*$ is the quantum parameter.  In [10] some  structural properties of $M_q(n)$ were obtained in a constructive-computational way.  For the usage of  Section 3, we quote some results of [10] in one proposition below. \v5

{\bf Proposition 1.3} [10, Corollary 2.2, Theorem 3.3, Theorem 4.2, Theorem 4.3(i)] Let $M_q(n)$ be the standard quantized matrix algebra in the sense of [2]. The following statements hold. \par
(1) The defining relations of $M_q(n)$ constitute a Gr\"obner-Shirshov basis and, consequently, $M_q(n)$ has a PBW $K$-basis which is in the standard form:
$$\B =\left\{\left. 1,~z^{k_{nn}}_{nn}z^{k_{nn-1}}_{nn-1}\cdots z^{k_{n1}}_{n_1}z^{k_{n-1n}}_{n-1n}\cdots z^{k_{n-11}}_{n-11}\cdots z^{k_{1n}}_{1n}\cdots z^{k_{11}}_{11}~\right |~k_{ij}\in \NZ,(i, j)\in I(n)\right\}.$$\par
(2) The Gelfand-Kirillov dimension GK.dim$M_q(n)=n^2$.\par
(3) The global homological dimension gl.dim$M_q(n)=n^2$.\par
(4) $M_q(n)$ is a solvable polynomial algebra in the sense of {[3], [7]), thereby is a Noetherian domain.\v5

\section*{2. $M_q(n)$ is an iterated skew polynomial algebra}
With notation and notions fixed in Section 1, this section is dedicated to show that $M_q(n)$ is an iterated skew polynomial K-algebra. \v5

For the sake of guaranteeing the argumentation below is convincible, we start by reviewing  briefly the definition of a skew polynomial ring here,  though this definition is well known. Recall that a skew polynomial ring is an associative ring $S$ over a ring $R$ in a variable $x$,  consisting of polynomials $f=\sum^n_{i=1}r_ix^i$ with coefficients in $R$, {\it where $x$ is not assumed to commute with the elements of} $R$, but is desired that{\parindent=1.4truecm\par
\item{(S1)} each polynomial should be expressed {\it uniquely} in the form $f=\sum_{i=1}^n r_ix^i$ for some $r_i\in R$ (or in other words, $S$ is a free (left) $R$-module with basis $1, x,x^2,\ldots ,x^k,\ldots$), and\par
\item{(S2)} $xr\in Rx+R$, i.e., $xr=\sigma (r)x+\delta (r)$ for some $\sigma(r)$, $\delta (r)\in R$.\par}{\parindent=0pt
It follows from (S1) and (S2) that $\sigma$ is a ring endomorphism of $R$ such that $\sigma (1)=1$, and that $\delta$ is a $\sigma$-derivation in the sense that $\delta\in$ End$_{\mathbb{Z}}R$ such that
$$\delta (r_1r_2)=\sigma (r_1)\delta (r_2)+\delta (r_1)r_2.$$
The skew polynomial ring $S$ as defined above is usually denoted by $S=R[x;\sigma ,\delta ]$.}\v5
To better understand the proof of Proposition 2.1 below, let us first look at how  $M_q(2)$ can be equipped with the structure of an iterated skew polynomial algebra. Note that $M_q(2)$ is defined subject to the relations
$$\begin{array}{ll} z_{11}z_{12}=qz_{12}z_{11},& z_{21}z_{22}=qz_{22}z_{21},\\
z_{11}z_{21}=qz_{21}z_{11},&z_{12}z_{22}=qz_{22}z_{12},\\
z_{12}z_{21}=z_{21}z_{12},&z_{11}z_{22}=z_{22}z_{11}+(q-q^{-1})z_{12}z_{21},\end{array}$$
and by Proposition 1.3, it has the PBW $K$-basis $\B =\{ z_{22}^{k_{22}}z_{21}^{k_{21}}z_{12}^{k_{12}}z_{11}^{k_{11}}~|~k_{ij}\in\NZ\}$.
Starting with the commutative polynomial ring $R_{22}=K[z_{22}]$, it is straightforward to verify that
$$M_q(2)=R_{22}[z_{21}; \sigma_{21}][z_{12}; \sigma_{12}][z_{11}; \sigma_{11},\delta_{11}],$$
where, with respect to the relations satisfied by the generators of $M_q(2)$, $\sigma_{21}$ is the algebra automorphism of $R_{22}=K[z_{22}]$ such that $\sigma_{21}(f(z_{22}))=f(qz_{22})$, $\sigma_{12}$ is the algebra automorphism of $K[z_{22}, z_{21}]$ such that$\sigma_{12}(f(z_{22}, z_{21}))=f(qz_{22}, z_{21})$, and $\sigma_{11}$ is the algebra  automorphism of $K[z_{22}, z_{21}, z_{12}]$ such that $\sigma_{11}(f(z_{22}, z_{21}, z_{12}))=f(z_{22}, qz_{21}, qz_{12})$, while $\delta_{11}$ is the $\sigma_{11}$-derivation of $K[z_{22}, z_{21}, z_{12}]$  such that $\delta_{11}(f(z_{22}, z_{21}, z_{12}))=f((q-q^{-1})z_{21}z_{12},  z_{21}, z_{12})$ (note that $z_{12}z_{21}=z_{21}z_{12}$).\v5

{\bf Proposition 2.1} The standard quantized matrix algebra $M_q(n)$ is an iterated skew polynomial ring of the form
$$\hskip -10truecm M_q(n)=K[z_{nn}, \ldots ,z_{n2},z_{n1}]$$
$$\hskip 2.1truecm [z_{n-1,n};\sigma^{(n)}_{n-1}][z_{n-1,n-1};\sigma^{(n-1)}_{n-1},\delta^{(n-1)}_{n-1}]
\cdots [z_{n-1,2};\sigma^{(2)}_{n-1},\delta^{(2)}_{n-1}]
[z_{n-1,1};\sigma^{(1)}_{n-1},\delta^{(1)}_{n-1}]$$
$$\hskip 2.1truecm [z_{n-2,n};\sigma_{n-2}^{(n)}][z_{n-2,n-1};\sigma_{n-2}^{(n-1)},\delta_{n-2}^{(n-1)}]
\cdots [z_{n-2,2};\sigma_{n-2}^{(2)}, \delta^{(2)}_{n-2}][z_{n-2,1};\sigma_{n-2}^{(1)},\delta_{n-2}^{(1)}]$$
$$\hskip -3.4truecm\cdots\hskip 1truecm\cdots\hskip 1truecm\cdots\hskip 1truecm\cdots\hskip 1truecm\cdots$$
$$\hskip -3.06truecm [z_{in};\sigma_i^{(n)}][z_{i,n-1};\sigma_i^{(n-1)},\delta_i^{(n-1)}]\cdots [z_{i1};\sigma_i^{(1)},\delta_i^{(1)}]$$
$$\hskip -3.4truecm\cdots\hskip 1truecm\cdots\hskip 1truecm\cdots\hskip 1truecm\cdots\hskip 1truecm\cdots$$
$$\hskip -3.06truecm [z_{1n};\sigma_1^{(n)}][z_{1,n-1};\sigma_1^{(n-1)},,\delta_1^{(n)}]\cdots [z_{11};\sigma_1^{(1)},\delta_1^{(1)}].$$
{\bf Proof} Bearing in mind Proposition 1.3 which says that $M_q(n)$ has the  PBW $K$-basis
$$\B =\left\{\left. 1,~z^{k_{nn}}_{nn}z^{k_{nn-1}}_{nn-1}\cdots z^{k_{n1}}_{n_1}z^{k_{n-1n}}_{n-1n}\cdots z^{k_{n-11}}_{n-11}\cdots z^{k_{1n}}_{1n}\cdots z^{k_{11}}_{11}~\right |~k_{ij}\in \NZ,(i, j)\in I(n)\right\},$$
and arranging the generators of $M_q(n)$ in the order
$$\begin{array}{rrrrr} z_{1n}&z_{1,n-1}&\cdots&z_{12}&z_{11}\\
\vdots~~&\vdots~~~&\cdots&\vdots~~&\vdots~~\\
z_{i-1,n}&z_{i-1,n-1}&\cdots&z_{i-1,2}&z_{i-1,1}\\
z_{in}&z_{i,n-1}&\cdots&z_{i2}&z_{i1}\\
z_{i+1,n}&z_{i+1,n-1}&\cdots&z_{i+1,2}&z_{i+1,1}\\
\vdots~~&\vdots~~~&\cdots&\vdots~~&\vdots~~\\
z_{n-1,n}&z_{n-1,n-1}&\cdots&z_{n-1,2}&z_{n-1,1}\\
z_{nn}&z_{n,n-1}&\cdots&z_{n2}&z_{n1}\end{array}$$
it follows from the relations $R_{[j<k]}$, $R_{[i<k]}$, $R_{[i<s, t<j]}$, and
$R_{[i<s, j<t]}$ that for $1\le i\le n$ we have
$$z_{i,n-p}z_{i,n-p+m}=qz_{i,n-p+m}z_{i,n-p}~\hbox{by}~R_{[j<k]}~\hbox{if}~1\le p\le n-1, 1\le m\le p,\eqno{(1)}$$
and for
$$\begin{array}{ll} i=n-1, n-2,\ldots ,1,& j=n,n-1,\ldots ,1,\\
k=1,2,\ldots ,n-1,& t=n,n-1,\ldots ,1,\end{array}$$
we have
$$z_{ij}z_{i+k,t}=\left\{\begin{array}{ll}
qz_{i+k,t}z_{ij},&\hbox{by}~R_{[i<k]}~\hbox{if}~t=j,\\
z_{i+k,t}z_{ij},&\hbox{by}~R_{[i<s,t<j]}~\hbox{if}~t<j,\\
z_{i+k,t}z_{ij}+(q-q^{-1})z_{it}z_{i+k,j}&\\
=z_{i+k,t}z_{ij}+(q-q^{-1})z_{i+k,j}z_{it},&\hbox{by}~R_{[i<s,j<t]}~\hbox{if}~j<t.
\end{array}\right.\eqno{(2)}$$
Thus, by means of  $(1)$ and (2) we may acquire the desired skew polynomial structure of $M_q(n)$ via a successive extension of skew polynomial subalgebras. More precisely, we first construct\par
$$A_n=K[z_{nn}][z_{n,n-1}; \sigma_n^{(n-1)}][z_{n,n-2};\sigma_n^{(n-2)}]\cdots [z_{n1};\sigma_n^{(1)}],$$
where for $1\le k\le n-1$, $\sigma_n^{(n-k)}$ is the algebra automorphism of $K[z_{nn}, z_{n,n-1},\ldots ,z_{n,n-k+1}]$ such that
$$\sigma_n^{(n-k)}(f(z_{nn}, z_{n,n-1},\ldots ,z_{n,n-k+1}))
=f(qz_{nn}, qz_{n,n-1},\ldots ,qz_{n,n-k+1}), $$
and then, we construct the skew polynomial subalgebra
$$A_{n-1}=A_n[z_{n-1,n};\sigma_{n-1}^{(n)}][z_{n-1,n-1};\sigma_{n-1}^{(n-1)},\delta_{n-1}^{(n-1)}]
[z_{n-1,n-2}: \sigma_{n-1}^{(n-2)},\delta_{n-1}^{(n-2)}]\cdots [z_{n-1,1};\sigma_{n-1}^{(1)},\delta_{n-1}^{(1)}]$$
where{\hskip 1.3truecm\par
\item{~} $\sigma_{n-1}^{(n)}$ is the algebra automorphism of $A_n=K[z_{nn},z_{n,n-1},\ldots ,z_{n1}]$, such that
$$\sigma_{n-1}^{(n)}(f(z_{nn},z_{n,n-1},\ldots ,z_{n1}))=f(qz_{nn}, z_{n,n-1},z_{n,n-2},\ldots , z_{n1});$$\par
\item{~} $\sigma_{n-1}^{(n-1)}$ is the algebra automorphism of $K[z_{nn}, \ldots ,z_{n1}, z_{n-1,n}]$, such that
$$\sigma_{n-1}^{(n-1)}(f(z_{nn},\ldots ,z_{n1}, z_{n-1,n}))=f(z_{nn}, qz_{n,n-1},z_{n,n-2},\ldots , z_{n1}, qz_{n-1,n}),$$
and $\delta_{n-1}^{(n-1)}$ is a $\sigma_{n-1}^{(n-1)}$-derivation of $K[z_{nn}, \ldots ,z_{n1}, z_{n-1,n}]$, such that
$$\delta_{n-1}^{(n-1)}(f(z_{nn},\ldots ,z_{n1}, z_{n-1,n}))=f((q-q^{-1})z_{n,n-1}z_{n-1,n}, z_{n,n-1},z_{n,n-2},\ldots , z_{n1}, z_{n-1,n}); $$\par
\item{~} $\sigma_{n-1}^{(n-2)}$ is the algebra automorphism of $K[z_{nn}, \ldots ,z_{n1}, z_{n-1,n}, z_{n-1,n-1}]$, such that
$$\begin{array}{l} \sigma_{n-1}^{(n-2)}(f(z_{nn},z_{n,n-1},\ldots ,z_{n1}, z_{n-1,n}, z_{n-1,n-1}))\\
=f(z_{nn}, z_{n,n-1},qz_{n,n-2},z_{n,n-3},\ldots , z_{n1}, qz_{n-1,n}, qz_{n-1,n-1}),
\end{array}$$
and $\delta_{n-1}^{(n-2)}$ is a $\sigma_{n-1}^{(n-2)}$-derivation of of $K[z_{nn}, \ldots ,z_{n1}, z_{n-1,n}, z_{n-1,n-1}]$  such that
$$\begin{array}{l} \delta_{n-1}^{(n-2)}(f(z_{nn},z_{n,n-1},\ldots ,z_{n1}, z_{n-1,n}, z_{n-1,n-1}))\\
=f((q-q^{-1})z_{n,n-2}z_{n-1,n}, (q-q^{-1})z_{n,n-2}z_{n-1,n-1}, z_{n,n-2}, z_{n,n-3},\ldots , z_{n1}, z_{n-1,n}, z_{n-1,n});\end{array}$$\par}
$$\hskip -3.4truecm\cdots\hskip 1truecm\cdots\hskip 1truecm\cdots\hskip 1truecm\cdots\hskip 1truecm\cdots$$
{\hskip 1.4truecm\par
\item{~} $\sigma_{n-1}^{(1)}$ is the algebra automorphism of $K[z_{nn}, \ldots ,z_{n1}, z_{n-1,n},\ldots ,z_{n-1},2]$, such that
$$\sigma_{n-1}^{(1)}(f(z_{nn},\ldots ,z_{n1}, z_{n-1,n},\ldots ,z_{n-1,2}))=f(z_{nn},\ldots , z_{n2}, qz_{n1},qz_{n-1,n},\ldots ,qz_{n-1,2}),$$
and $\delta_{n-1}^{(1)}$ is a $\sigma_{n-1}^{(1)}$-derivation of $K[z_{nn}, \ldots ,z_{n1}, z_{n-1,n},\ldots ,z_{n-1},2]$, such that
$$\begin{array}{l}\delta_{n-1}^{(1)}(f(z_{nn},\ldots ,z_{n1}, z_{n-1,n},\ldots ,z_{n-1,2}))\\
=f((q-q^{-1})z_{n1}z_{n-1,n},\ldots , (q-q^{-1})z_{n1}z_{n-1,2} ,z_{n1},\ldots , z_{n-1,2}). \end{array}$$\par}{\parindent=0pt
In general suppose that the skew polynomial structure of the subalgebra
$$A_{i+1}=K[z_{nn},\ldots ,z_{n1}, z_{n-1,n},\ldots ,z_{n-1,1}\ldots ,z_{i+1,n},\ldots ,z_{i+1,1}]$$
is already established. We then may construct the skew polynomial subalgebra
$$A_i=A_{i+1}[z_{in};\sigma_i^{(n)},\delta_i^{(n)}][z_{i,n-1};\sigma_i^{(n-1)}]\cdots [z_{i1};\sigma_i^{(1)},\delta_i^{(1)}],$$
where}{\hskip 1.4truecm\par
\item{~} $\begin{array}{l} \sigma_i^{(n)}(f(z_{nn},\ldots ,z_{n1},z_{n-1,n},\ldots ,z_{n-1,1},\ldots , z_{i+1,n},\ldots ,z_{i+1,1}))\\
    =f(qz_{nn}, z_{n,n-1},\ldots ,z_{n1}, qz_{n-1,n}, z_{n-1,n-1},\ldots ,z_{n-1,1},\ldots , qz_{i+1,n}, z_{i+1, n-1},\ldots ,z_{i+1,1}),\end{array}$\par
\item{~} $\begin{array}{l} \delta_i^{(n)}(f(z_{nn},\ldots ,z_{n1},z_{n-1,n},\ldots ,z_{n-1,1},\ldots , z_{i+1,n},\ldots ,z_{i+1,1}))\\
    =f((q-q^{-1})z_{n,n-1}z_{in},, z_{n,n-1},\ldots ,z_{n1}, (q-q^{-1})z_{n-1,n-1}z_{in}, z_{n-1,n-1},\ldots ,z_{n-1,1},\\
    \quad~~\ldots , (q-q^{-1})z_{i+1,n-1}z_{in}, z_{i+1, n-1},\ldots ,z_{i+1,1}, z_{in});\end{array}$
\item{~} $\begin{array}{l} \sigma_i^{(n-1)}(f(z_{nn},\ldots ,z_{n1},z_{n-1,n},\ldots ,z_{n-1,1},\ldots , z_{i+1,n},\ldots ,z_{i+1,1}, z_{in}))\\
    =f(z_{nn}, qz_{n,n-1},z_{n,n-2}\ldots ,z_{n1}, z_{n-1,n}, qz_{n-1,n-1},z_{n-1,n-2}\ldots ,z_{n-1,1},\\
    \quad~~\ldots , z_{i+1,n}, qz_{i+1, n-1},z_{i+1,n-2}\ldots ,z_{i+1,1},qz_{in}),\end{array}$,
\item{~} $\begin{array}{l} \delta_i^{(n-1)}(f(z_{nn},\ldots ,z_{n1},z_{n-1,n},\ldots ,z_{n-1,1},\ldots , z_{i+1,n},\ldots ,z_{i+1,1}), z_{in})\\
    =f((q-q^{-1})z_{n,n-1}z_{in},, z_{n,n-1}, z_{n,n-2},\ldots ,z_{n1}, (q-q^{-1})z_{n-1,n-1}z_{in}, z_{n-1,n-1},z_{n-1,n-2},\ldots ,z_{n-1,1},\\
    \quad~~\ldots , (q-q^{-1})z_{i+1,n-1}z_{in}, z_{i+1, n-1},z_{i+1,n-2},\ldots ,z_{i+1,1}, z_{in});\end{array}$ \par}
$$\hskip -3.4truecm\cdots\hskip 1truecm\cdots\hskip 1truecm\cdots\hskip 1truecm\cdots\hskip 1truecm\cdots$$
{\parindent=.6truecm\par
\item{~} $\begin{array}{l} \sigma_i^{(1)}(f(z_{nn},\ldots ,z_{n2},z_{n1},\ldots ,z_{n-1,n}, \ldots z_{n-1,2},z_{n-1,1},\ldots , z_{i+1,n},\ldots ,z_{i+1,2},z_{i+1,1}, z_{in},\ldots ,z_{i3}, z_{i2}))\\
    =f(z_{nn},\ldots ,z_{n2}, qz_{n1},z_{n-1,n}, \ldots ,z_{n-1,2}, qz_{n-1,1},\ldots ,z_{i+1,n},\ldots , z_{i+1,2}, qz_{i+1,1},
    qz_{in},\ldots , qz_{i3}, qz_{i2}),\end{array}$
\item{~} $\begin{array}{l} \delta_i^{(1)}(f(z_{nn},\ldots ,z_{n2},z_{n1},\ldots ,z_{n-1,n}, \ldots z_{n-1,2},z_{n-1,1},\ldots , z_{i+1,n},\ldots ,z_{i+1,2},z_{i+1,1}, z_{in},\ldots ,z_{i3}, z_{i2}))\\
    =f((q-q^{-1})z_{n1}z_{in}, (q-q^{-1})z_{n1}z_{i,n-1}, z_{n,n-2},\ldots ,z_{n1},\\
    ~~~~~~(q-q^{-1})z_{n-1,1}z_{in}, (q-q^{-1})z_{n-1,1}z_{i,n-1}, z_{n-1,n-2},\ldots ,z_{n-1,1},\\
    ~~~~\ldots (q-q^{-1})z_{i+1,1}z_{in}, (q-q^{-1})z_{i+1,1}z_{i,n-1},z_{i, n-2},\ldots ,z_{i2}).
    \end{array}$\v5}{\parindent=0pt
So, repeating the above iterating procedure for a finite number of times, the desired skew polynomial structure of $M_q(n)$ is finally obtained.\QED}\v5

\section*{3. Main result}
With notation and notions as before, the aim of this section is to prove the main result of this paper, which is mentioned below.\v5

{\bf Theorem 3.1} Let $M_q(n)$ be the standard quantized matrix algebra. The following statements hold.\par
(i) $M_q(n)$ is an Auslander regular algebra.\par
(ii) $M_q(n)$ is Cohen-Macaulay.\par
(iii) $M_q(n)$ is an Artin-Schelter regular algebra.\par
(iv) $M_q(n)$ is a maximal order in its quotient division algebra.\v5

To this end, let us first recall from [8] a basic result concerning the Auslander regularity of an iterated skew polynomial ring.\v5

{\bf Proposition 3.2} [8, P.170, Theorem 6(1), Theorem 6(4)] Suppose that a ring $R$ is left and right Noetherian and Auslander regular, then the following two statements hold.\par
(i) Any skew polynomial ring $R[t;\sigma]$ is left and right Noetherian and Auslander regular, where $\sigma$ is an automorphism of $R$.\par
(ii) Any skew polynomial ring $R[t; \sigma , \delta]$  is left and right Noetherian and Auslander regular, where $\sigma$ is an automorphism of $R$, and $\delta$ is a $\sigma$-derivation of $R$. \v5

{\bf Proof of Theorem 3.1} (i) Since $M_q(n)$ is left and right Noetherian with finite homological dimension $n^2$ by Proposition 1.3, and by the proof of Proposition 2.1, the subalgebra $A_n=K[z_{nn}, z_{n,n-1}, \ldots , z_{n1}]$ of $M_q(n)$ is an iterated skew polynomial algebra over the commutative polynomial algebra $K[z_{nn}]$ via a series of algebra automorphisms $\sigma_n^{(k)}$, $1\le k\le n-1$,  while $M_q(n)$ is an iterated skew polynomial algebra over $A_n$ via a series of algebra automorphism $\sigma_i^{(k)}$ and $\sigma_i^{(k)}$-derivations, $1\le k\le n$, it follows from Proposition 3.2 that $A_q(n)$ is an Auslander regular algebra.\par
(ii) Due to the fact that $M_q(n)$ is left and right Noetherian with gl.dim$M_q(n)=n^2=$ GK.dim$M_q(n)$ (Proposition 1.3), and it is an iterated skew polynomial algebra (Proposition 2.1), it follows from [8, CH.III] and [4] that $M_q(n)$ is Cohen-Macaulay.\par
(iii) Note that $M_q(n)$ is obviously a connected $\NZ$-graded $K$-algebra in which every generator has degree 1. Since $M_q(n)$ is a (Noetherian) Auslander regular algebra by (i), this assertion may follow from [4] (see also a remark given in [11, PP. 127--128]).\par
(iv) First note that $M_q(n)$ is a solvable polynomial algebra, thereby it is a Noetherian domain by Proposition 1.3, and it is also stably free by [7, Theorem 3.3.3]. Moreover,  since $M_q(n)$ is Auslander regular and Cohen-Macaulay by assertions (i) and (ii),  it follows from Proposition 1.2 that $M_q(n)$ is  a maximal order in its quotient division algebra $Q(A)$.\par
Summing up, the theorem is proved.\QED \v5

{\bf Acknowledgement}\par
The author is very grateful to professor Huishi Li for proposing the research topic, for very valuable discussions, and for his patient and meticulous guidance during preparing the manuscript. \v5

\centerline{Refeerence}{\parindent=1.47truecm\par

\item{[1]} Artin, M., Schelter, W. (1987). Graded algebras of global dimension 3. {\it Adv. Math}. 66:171--216.

\item{[2]} Faddeev, L. D.,   Reshetikhin, N. Yu., Takhtajan, L. A. (1988).  Quantization of Lie groups and Lie algebras. {\it Algebraic Analysis}, Academic Press, 129-140.\par

\item{[3]} Kandri-Rody A., Weispfenning, V. (1990).  Non-commutative
Gr\"obner bases in algebras of solvable type. {\it J. Symbolic
Comput.}, 9:1--26. Also available as: Technical Report University of Passau,
MIP-8807, March 1988.

\item{[4]} Levasseur, T. (1992). Some properties of non-commutative regular rings. {\it Glasgow Math. J}. 34:277--300.

\item{[5]} Li, H. (2011). {\it Gr\"obner Bases in Ring Theory}. World Scientific Publishing Co. \url{https://doi.org/10.1142/8223}\par

\item{[6]} Li, H. (2018). An elimination lemma for algebras with PBW bases. {\it Communications in
Algebra}, 46(8):3520-3532.\par

\item{[7]} Li, H. (2021). {\it Noncommutative polynomial algebras of solvable type and their modules: Basic constructive-computational theory and methods}. Chapman and Hall/CRC Press. \url{https://doi.org/10.1201/9781003213192}\par

\item{[8]} Li, H., Van Oystaeyen, F. (1996, 2003). {\it Zariskian Filtrations}. K-Monograph in Mathematics, Vol.2. Kluwer Academic Publishers, Berlin Heidelberg: Springer-Verlag.
    \url{https://doi.org/10.1007/978-94-015-8759-4_2}\par

\item{[9]} Stafford, J. T. (1994).  Auslander-regular algebras and maximal orders.  {\it J. London Math. Soc}. (2)50:276--292.\par

\item{[10]} Tuniyaz, R. (2023). On the standard quantized
matrix algebra $M_q(n)$: from a constructive-computational viewpoint.
{\it Journal of Algebra and Its Applications}. 18 Jan. 2023 Published
online: https://doi.org/10.1142/S0219498824500920.

\item{[11]} Van Oystaeyen, F. (2000). {\it Algebraic Geometry for Associative Algebras}. (Chapman \& Hall/CRC Pure and Applied Mathematics) CRC Press.

\end{document}